\documentclass[greek,english,12pt,reqno]{amsart}
\usepackage{amssymb}
\usepackage{pstricks}

\newtheorem{theorem}{Theorem}[section]
\newtheorem{lemma}[theorem]{Lemma}

\newtheorem{proposition}[theorem]{Proposition}

\theoremstyle{definition}

\newtheorem{definition}[theorem]{Definition}

\newtheorem{example}[theorem]{Example}

\newtheorem{remark}[theorem]{Remark}

\def\ord{\operatorname{ord}}

\def\Z{\mathbb{Z}}

\def\Z_S{{\mathbb Z_S}}
\def\beq {\begin{equation}}
\def\endeq {\end{equation}}

\def\eds{elliptic divisibility sequence}
\def\edssp{elliptic divisibility sequence }

\def\edsssp{elliptic divisibility sequences }

\begin{document}
\title{Descent on elliptic curves and Hilbert's Tenth Problem}
\author{Kirsten Eisentr\"ager}
\address{Department of Mathematics, The Pennsylvania State University,
University Park, PA 16802, USA.}
\email{eisentra@math.psu.edu}
\author{Graham Everest}
\address{School of Mathematics, University of East Anglia,
Norwich NR4 7TJ, UK.}\email{g.everest@uea.ac.uk}
\thanks{The authors thank the ICMS in Edinburgh for the workshop on Number Theory and
Computability in 2007, funded by EPSRC and the LMS.}
\thanks{The first author
was partially supported by NSF grant DMS-0801123 and a grant from the
John Templeton Foundation.}


\begin{abstract}
Descent via an isogeny on an elliptic curve is used to
construct two subrings of the field of rational numbers, which are
complementary in a strong sense, and
for which Hilbert's Tenth Problem is
undecidable. This method further develops that of Poonen, who
used elliptic divisibility sequences to obtain
undecidability results for some large subrings of the rational numbers.
\end{abstract}
\maketitle

\section{Hilbert's Tenth Problem}

In 1970, Matijasevi\v{c} \cite{ym}, building upon earlier work of
Davis, Putnam and Robinson \cite{dpr}, resolved negatively
Hilbert's Tenth Problem for the ring~$\mathbb Z$, of rational
integers. This means there is no general algorithm which will
decide if a polynomial equation, in several variables, with integer
coefficients has an integral solution. Equivalently, one says
Hilbert's Tenth Problem is {\it undecidable}
for the integers.
See~\cite[Chapter 1]{as} for a full overview and background
reading. The same problem, except now over the rational
field~$\mathbb Q$, has not been resolved. In other words, it is
not known if there is an algorithm which will decide if a
polynomial equation with integer coefficients (or
rational coefficients, it doesn't matter) has a rational solution.

Recently, Poonen \cite{bjorn} took a giant leap in this direction
by proving the same negative result for some large subrings
of~$\mathbb Q$. To make this precise, given a prime $p$ of
$\mathbb Z$, let $|.|_p$ denote the usual $p$-adic absolute value.
Let $S$ denote a set of rational primes. Write
$$\mathbb Z_S=\mathbb Z[1/S]=\{x\in \mathbb Q:|x|_p\leq 1 \mbox{
for all }p\notin S\},
$$
for the ring of $S$-integers of $\mathbb Q$.
\begin{theorem}[Poonen \cite{bjorn}]\label{PT} There
are recursive sets $S$ of primes having density 1 with the property that Hilbert's Tenth
Problem for~$\mathbb Z_S$ is undecidable.
\end{theorem}

Given the importance of Theorem\ref{PT}, it is surely worth
investigating more closely the subrings~$\mathbb Z_S$ of~$\mathbb
Q$ to which it applies. Besides the
intrinsic interest, the hope remains that a solution for the
rational field might be accessed through the rings~$\mathbb Z_S$. The aim of this
paper is to extend Poonen's method in a non-trivial fashion, by using
descent on elliptic curves. As motivation, note that although the rings~$\mathbb Z_S$ in
Theorem~\ref{PT} are formed by inverting sets
of primes with density~$1$, the sets of primes are necessarily co-infinite.
It is not even clear from \cite{bjorn}
whether a finite collection of such rings will generate~$\mathbb
Q$. In this paper we provide examples where two rings
suffice, and which are complementary in a strong sense. Write $\mathbb P$ for the set of all prime numbers.
\begin{definition} Two subsets of $\mathbb P$ are said to be {\it complementary} if their
union is~$\mathbb P$. They are said to be {\it exactly complementary} if, in addition, they
have empty intersection. A subset of~$\mathbb Z$ is said to be {\it recursive}
if there is an algorithm to decide if any given integer lies in that subset.
\end{definition}

\begin{theorem}\label{h10} There are exactly complementary recursive sets $S,T\subset \mathbb P$
such that Hilbert's Tenth Problem is undecidable for both
rings~$\mathbb Z_{S}$ and~$\mathbb Z_T$.
\end{theorem}
Given sets $S$ and $T$ as in Theorem~\ref{h10},
any element $q\in \mathbb Q^*$ can be written
\begin{equation}\label{da}
q=st \mbox{ with } s\in \mathbb Z_S^*, t \in \mathbb Z_T^*,
\end{equation}
in a way that is unique up to sign.
Equation (\ref{da}) is a kind of {\it diophantine definition} (see \cite[Chapter 1]{as}) of the product group
$$\mathbb Z_S^* \times \mathbb Z_T^*
$$
over the group $\mathbb Q^*$. Hitherto, the concept of diophantine definition has
only been studied for rings and it is not known whether the property
in (\ref{da}) will
permit some kind of `lifting' of undecidability to the
rational field.

Theorem \ref{PT} was proved by
constructing a {\em diophantine model} of the positive integers in the
ring~$\mathbb Z_S$ using integer sequences ({\it elliptic divisibility sequences})
constructed from elliptic curves.
Consult \cite{bjorn}
and \cite[Chapter 12]{as} for the background and full details of the definitions.

\begin{definition}\label{model}
We say a set $A\subseteq \mathbb Z_S$ is {\it diophantine} over $\mathbb Z_S$ to mean $A$ is a projection
of the set of solutions of a diophantine
equation over~$\mathbb Z_S$. In other words, there exists $f\in \mathbb Z_S[y,x_1,\dots ,x_n]$ such that
$$a\in A \iff \exists \underline t \in \mathbb Z_S^n \mbox{ with }
f(a,t_1,\dots ,t_n)=0.
$$
A {\it diophantine model} of $\mathbb N$ over $\mathbb Z_S$ is a bijection
$\mathbb N \longleftrightarrow A$, where $A$ is diophantine over
$\mathbb Z_S$ with the additional property that the graphs of $+$ and~$\times$
correspond to diophantine subsets of $A^3$.
\end{definition}

As in \cite{bjorn,as} the undecidability of $\mathbb Z$ is essentially
equivalent to that of~$\mathbb N$, together with~$+$ and~$\times$.
Technically it is easier to
model this latter set.
Definition~\ref{model} is important because it allows undecidability results to be lifted
from the positive integers to the set~$A$. Exactly the same issue arises in this paper. On this point
the arguments are identical.

\subsection{Background Results}
Definition \ref{model} brings to the fore the role played by diophantine equations.
What follows is a brief overview of earlier results, which shows how
advances have been made by changing the underlying equation.
In~\cite{kr}, Kim and Roush resolved Hilbert's Tenth Problem negatively for
rings~$\mathbb Z_S$ when~$S$ consists of a single prime. The underlying
equation involved a
quadratic form, much in the spirit of earlier work by Julia
Robinson. A negative answer to Hilbert's Tenth Problem for rings~$\mathbb Z_S$, when~$S$ is
finite, follows using the concept of {\it diophantine class}
as in~\cite[Chapter 4]{as}.
Shlapentokh \cite{asimrn} resolved
Hilbert's Tenth Problem problem for some
large subrings of number fields, where the underlying diophantine equation arose from a
homogeneous polynomial known as a {\it norm form}.
Poonen's Theorem~\cite{bjorn}, Theorem~\ref{PT}, is important because, for the first time, it resolved
Hilbert's Tenth Problem for certain rings~$\mathbb Z_S$ when~$S$
is infinite. His underlying equation was an elliptic curve and the
definition of diophantine model was satisfied by using an \eds. Theorem~\ref{PT}
has been generalized
to subrings of number fields in \cite{poonenshalpentokh}.
In another interesting direction, Cornelissen and
Zahidi \cite{guntherkarim} also
used an \edssp to obtain decidability results.

\section{Elliptic Curves}

Let $E$ denote an elliptic curve,
\begin{equation}\label{defE}
E: \quad y^2+a_1xy+a_3y=x^3+a_2x^2+a_4x+a_6
\end{equation}
where $a_1,\dots ,a_6$ denote integers.  Consult
\cite{cassels,silbook1} for the basic properties of elliptic
curves. Suppose $Q\in E(\mathbb Q)$ denotes a non-torsion rational
point. The shape of the defining equation (\ref{defE}) forces the
denominator of the $x$-coordinate of a rational point to be a square, and that
of the $y$-coordinate to be a cube.
For $1\leq n \in \mathbb N$, write $nQ$ for the $n$-th multiple of~$Q$ according to the usual
addition law on~$E$. Then
\begin{equation}\label{defEDS}
nQ=\left(\frac{A_n}{B_n^2},\frac{C_n}{B_n^3}\right),
\end{equation}
with $A_n,B_n,$ and $C_n$ denoting integers which satisfy $B_n>0$ and
gcd$(B_n,A_nC_n)=1$.
The sequence~$B=(B_n)$ is known as an {\em elliptic divisibility
sequence}. An important property of $B$ (the `divisibility'
part of its name) is the following
\begin{equation}\label{divseq}
n|m \mbox{ in $\mathbb N$ implies } B_n|B_m \mbox{ in $\mathbb
N$}.
\end{equation}

\begin{definition}\label{pd}
Let $B=(B_n)$ denote a sequence with integer terms. We say an
integer $d>1$ is a {\it primitive divisor} of the term $B_n\neq 0$ if
\begin{enumerate}\item $d\mid B_n$ and \item gcd$(d,B_m)=1$ for all non-zero
terms~$B_m$ with $0<m<n$.\end{enumerate}\end{definition}

In 1986, Silverman \cite{silprimdiv} proved an analogue of Bang's
theorem \cite{bang}, that the terms of \edsssp have primitive divisors
for all sufficiently large indices.

\begin{remark}\label{pdasorder}If $l$ denotes any prime divisor of $d$
as in Definition~\ref{pd} then it is referred to a {\it primitive
prime divisor} of~$B_n$. Provided~$l$ is a prime of non-singular
(hereafter {\it good}) reduction for~$E$, an important group theoretic interpretation of
the situation is that~$n$ is the order of the point~$Q$ mod~$l$ on
the reduced curve. It follows that
\begin{equation}\label{divrel}B_m \equiv 0 \mbox{ mod } l \mbox{ if and only if } n|m.
\end{equation}
\end{remark}

\subsection{Two Primitive Divisors}
Silverman's Theorem ensures that for all sufficiently large~$n$,
every term~$B_n$ has a primitive divisor. More can be said when descent via an
isogeny is possible.

\begin{proposition}\label{foranyt}Suppose $Q\in E(\mathbb Q)$ is a non-torsion
point which generates an \edssp $B=(B_n)$. Let~$\sigma:E'\to E$
denote an isogeny of prime degree~$q$ and assume
~$Q=\sigma(Q')$ for some rational point~$Q'\in E'(\mathbb Q)$.
Then all terms~$B_n$, with~$n$ sufficiently large
and coprime to~$q$, have at least\footnote{`At least' means no
fewer than.}~$2$ distinct primitive prime divisors.
\end{proposition}

The techniques needed
for the proof of Proposition~\ref{foranyt} draw upon those used in \cite{primeds}, \cite{pe} and \cite{emw}.
Write $\sigma^*:E\to E'$ for the dual isogeny, also of prime degree~$q$. The
compositions $\sigma \sigma^*$ and $\sigma^* \sigma$ are the maps $[q]$ (times $q$) on
$E$ and $E'$ respectively.
Given $\sigma:E'\to E$
and~$Q=\sigma(Q')$, write~$b=(b_n)$ for the \edssp corresponding to~$Q'$.

\begin{definition} Given any \edssp $B=(B_n)$, write $B_n^*$ for the {\it primitive part}
of $B_n$. This is the maximal divisor of~$B_n$ which is coprime to all the terms~$B_m$
with $0<m<n$.\end{definition}

\begin{lemma}\label{rates} There are positive constants $h$ and $h'$ with $h=qh'$ such that for large~$n$:
\begin{itemize}\item $\log B_n \sim hn^2$ \item $\log b_n \sim h'n^2$ \item $\log B_n^* \geq .547hn^2$
\end{itemize}
\end{lemma}
\begin{proof}[\sc Proof of Lemma \ref{rates}] The first two formulae represent a strong form
of Siegel's Theorem \cite[Chapter IX]{silbook1}. See also \cite{pe} for a direct proof using elliptic
transcendence theory. For the third formula, use \cite[Lemma 3.3]{emw}. From this follows
a lower bound of the form
$$\log B_n^* \geq \log B_n - \sum_{p|n}\log \left(p^2B_{\frac{n}{p}}\right).
$$
Now apply the first formula together with the upper bound
$$\sum_{p}\frac{1}{p^2}<.453
$$
The constants $h=\hat h(Q)$ and $h'=\hat h(Q')$ represent the canonical height of $Q$ and $Q'$ on their
respective curves. The relation $h=qh'$ is a property of the canonical height under
isogeny.
\end{proof}

\begin{proof}[\sc Proof of Proposition \ref{foranyt}]  Let~$p$
denote any prime of non-singular reduction for~$E$ (or~$E'$, the two curves
share the same set of good reduction primes). By applying
an isomorphism to $E'$ if necessary, Velu's formulae \cite{velu} imply
\begin{equation}\label{divdiv}
\ord_p(b_n)\leq \ord_p(B_n)\leq \ord_p(b_{qn}) \mbox{ for all } n \geq 0.
\end{equation}
For all sufficiently large $n$, the term $b_n$ has a primitive prime divisor~$l_n$.
Assume that~$n$ is
large enough to guarantee that~$l_n$ is a prime of good
reduction.
Then~$l_n$ is a divisor of $B_n$ by (\ref{divdiv}).
If gcd$(q,n)=1$, we claim~$l_n$ is actually a primitive
prime divisor of~$B_n$. If not, then
$l_n|B_m,$ for some $0<m<n$, chosen minimally. In group-theoretic terms, see Remark~\ref{pdasorder}, this
means $$mQ \equiv nQ' \equiv O \mbox{ mod } l_n$$
on the corresponding reduced curves. Now (\ref{divrel}) implies the following
divisibility relations
$$m|n \mbox{ and } n|qm.
$$
Since $q$ is prime these force $n=qm$, contradicting the assumption that gcd$(q,n)=1$.

The proof is now completed using the data about growth rates of the various
sequences in Lemma~\ref{rates}. Note \cite{primeds,pe} that the contribution to~$B_n$ from
primes of singular (hereafter {\it bad}) reduction is negligible. In particular, it follows
that $B_n^*$ grows asymptotically
faster than its divisor~$b_n^*$. It remains to show that for all large enough~$n$, $B_n^*$ has a prime
divisor which is coprime to~$b_n^*$. This forces $B_n^*$ to have at least two distinct
prime divisors, which is the desired conclusion. To prove this last claim use the
following property of \edsssp from the $p$-adic
theory of elliptic curves. The property can be sourced in \cite[Chapter IV]{silbook1} and is stated
in \cite[Lemma 3.1]{emw}. It says that provided $l_n>2$ and $l_n|b_n$,
\begin{equation}\label{ordord}
\ord_{l_n}(b_{qn})=\ord_{l_n}(b_n)+\ord_{l_n}(q).
\end{equation}
If $l_n|\gcd(B_n^*,B_n^*/b_n^*)$ then (\ref{divdiv}) implies that
$$\ord_{l_n}(b_{qn})>\ord_{l_n}(b_n).
$$
Now (\ref{ordord}) shows that the only way this can happen is if $l_n|q$. We may assume
$n$ is large enough to avoid this possibility.
\end{proof}

\section{Proof of Theorem \ref{h10}}

\begin{proof}[\sc Proof of Theorem \ref{h10}]
Assume $E$ is an elliptic curve and there is an isogeny $\sigma: E'\to E$ of
prime degree~$q$ such that:
\begin{itemize}
\item $E(\mathbb Q)=<Q>\simeq \mathbb Z,$
\item $E(\mathbb R)$ has only one real-connected
component,
\item $Q$ is the image of a $\mathbb Q$-rational point under~$\sigma$,
\item $B_q>1$.
\end{itemize}
\begin{example}\label{exists}$y^2=x^3-4 \quad Q=[2,2]$
This curve has conductor $432$ and it appears as $b1$ in Cremona's tables \cite{jec}.
There is a $3$-isogeny\footnote{An example meeting the needs of the proof
will necessarily have degree greater
than 2. This is because a curve with a 2-isogeny will have a
rational 2-torsion point.}
from the curve (called $b2$) $y^2=x^3+108$
which maps $[6,18]$ to $Q$. The properties claimed are easily checked.
\end{example}

By Proposition \ref{foranyt}, for all
sufficiently large primes~$l$, $B_l$ has at least two distinct primitive
prime divisors. (Any prime divisor of a term~$B_l$, with~$l$ a
prime, is necessarily a primitive prime divisor, using
(\ref{divrel}): the essential contribution of
Proposition~\ref{foranyt} is that it guarantees at least two distinct prime
divisors). Also, by Proposition \ref{foranyt}, each
term~$B_{ll'}$, where $l,l'$ are distinct primes, has at least two
distinct primitive prime divisors, except possibly for a finite number of
pairs~$(l,l')$, provided the primes~$l$ and~$l'$
are distinct from~$q$.
\begin{definition}\label{ab}For every prime~$l$, let~$a_l\geq 1$ denote the
smallest integer such that~$B_{l^{a_l}}>1$.
Let $L$ denote the set of primes~$l$ such that~$a_l>1$. Then
$L$ is finite by Siegel's Theorem. In Example~\ref{exists}, $2\in L$
because $B_2=1$. Also~\cite[Proposition 2.5]{silbook1}, the
set of everywhere good reduction points form a subgroup of~$E(\mathbb Q)$. It follows that for
each bad reduction prime~$p$, this
subgroup contains the kernel of reduction mod~$p$.
Thus~$b_p>1$ exists such that $p|B_n$ if and only if $b_p|n$. Write~$b$
for the computable number consisting of the largest of the~$b_p$.
\end{definition}

Exactly as in \cite{bjorn}, use Vinogradov's Theorem \cite[Chapter
XI]{vin} on the additive circle $$E(\mathbb R)\simeq \mathbb
R/\mathbb Z\simeq [0,1).$$

This theorem guarantees that the
multiples~$lQ$, with~$l$ prime, are dense in the real
curve~$E(\mathbb R)$.
With $b$ as in Definition \ref{ab}, choose a set~$U$ of primes
inductively as follows: given $l_1,\dots ,l_{i-1}$ choose $l_i$ to be the
smallest prime outside $L\cup \{q\}$ with $l_i>l_j>b$ for all $j<i$ and
\begin{equation}\label{ineq}
|y_i-i|<1/10i, \mbox{ where $l_iQ=(x_i,y_i)$.}
\end{equation}

\subsection{Definition of $S$} For all sufficiently large $n$ define
$p_n$ to be the largest primitive prime divisor of~$B_n$, which is also
a good reduction prime: this exists by~\cite{silprimdiv}. Let
the set~$S_1$ consist of the
prime divisors of all
the terms~$B_{l_i}, i\in\mathbb N$. The elements of~$S_1$ are all good reduction primes
because $l_i>b$ for all $i\in \mathbb N$.
Now define the set~$S_2$ by
\begin{equation}\label{S'}S_2=\{p_l:l \mbox{ prime}\neq l_i \forall
i\}\cup\{p_{l_il_j}:1\leq j\leq i\}\cup \{p_{ll_i}:l\in L, i\in \mathbb N\}.
\end{equation}
Clearly $S_1\cap S_2=\emptyset$. Let $S$ denote any set
containing~$S_1$ but disjoint from~$S_2$.
The primes in~$S_2$ act as witnesses to elements being outside
of~$E(\mathbb Z_S)$. In other words, as in \cite{bjorn,as},
\begin{equation}\label{ohyes}
\cup_i \{\pm l_iQ\}=E(\mathbb Z_S),
\end{equation}
with at most finitely many exceptions. For convenience, a proof of
(\ref{ohyes}) is now
indicated. Clearly all $\pm l_iQ\in E(\mathbb Z_S)$ because the
primes dividing terms~$B_{l_i}$ lie in~$S_1 \subseteq S$.
On the other hand, for all large enough~$n$,

\begin{eqnarray*}
n \neq \pm l_i \mbox{ some } i & \implies &
l|n \mbox{ for some } l \neq l_i \mbox{ or } l_il_j | n \mbox{ or } ll_i|n,
l\in L \\ &\implies & \exists p|B_n \mbox{ with } p=p_l \mbox{ or }
p=p_{l_il_j}\mbox{ or } p=p_{ll_i}, l\in L\\
& \implies & \exists p|B_n
\mbox{ with } p\in S_2 \subseteq S'\\
& \implies & nQ \notin E(\mathbb Z_S).
\end{eqnarray*}

Write $A$ for the set $A=\{y_i:l_iQ=(x_i,y_i)\}$.
The bijection required by Definition \ref{model}
is $i\leftrightarrow y_i$. Plainly $A$ is diophantine over~$\mathbb Z_S$, using the underlying
diophantine equation of the elliptic curve.
\begin{lemma}\label{plustimes}The graphs of $+$ and $\times$ correspond to diophantine subsets of~$A^3$.
\end{lemma}
Lemma~\ref{plustimes} is proved in \cite[Section 10]{bjorn}.
For example, it follows from~(\ref{ineq}) that $m+n-q$ differs from $y_m+y_n-y_q$
by at most 3/10. Therefore adding on~$\mathbb N$ corresponds to adding
on~$A$ then rounding to the nearest element. In other words, $m+n=q$ corresponds to a diophantine
predicate on~$A$. Multiplication is similar because it can be obtained by squaring
and adding. It follows {\it mutatis
mutandis} that~$\mathbb N$ has a diophantine model in~$\mathbb
Z_S$ and therefore Hilbert's Tenth Problem is undecidable
in~$\mathbb Z_S$.

\subsection{Definition of $T$} Define $p_n'$, for all sufficiently
large $n$ with~$q\nmid n$, to
be the second largest good reduction primitive prime divisor of~$B_n$. This exists by
Proposition~\ref{foranyt}. Now define
$T_1=S_1$ and
\begin{equation}\label{T'}T_2=\{p_l':l_i \neq l \mbox{ prime} \}
\cup\{p_{l_il_j}':1\leq j\leq i\}\cup\{p_{ll_i}':l\in L, i\in \mathbb N\}.
\end{equation}
The hypothesis that $B_q>1$ implies $q\notin L$. This is used
to guarantee that $p_{ll_i}'$ exists for all large $i$.
Let~$T$ denote any set containing~$T_1$ but disjoint from~$T_2$.
In exactly
the same way as before, $$\cup_i \{\pm l_iQ\}=E(\mathbb Z_T),$$
with at most finitely many exceptions. Again~$\mathbb N$ has a
diophantine model in~$\mathbb Z_T$ and therefore Hilbert's Tenth Problem is undecidable
in~$\mathbb Z_T$. Note that $S_2\cap T_2=\emptyset$ so choose
$$S=\mathbb P-S_2 \mbox{ and } T=\mathbb P-T_2.
$$
This results in $S\cup T=\mathbb P$. Subsequently, it will be argued that
the sets~$S_i,T_i,i=1,2$ are recursive. It follows that both~$S$ and~$T$
can be chosen to be recursive. This completes the proof that complementary
sets can be found. This argument will now be refined.

\subsection{Exactly Complementary Sets} To show that $S$ and $T$ may be chosen in an exactly
complementary fashion, choose
the sets $S_1$ and $S_2$ exactly as before.
Now choose a set~$U'$ of primes
inductively as follows: given $l_1',\dots ,l_{i-1}'$ choose $l_i'$ to be the
smallest prime outside $U\cup L\cup \{q\}$ with $l_i'>l_j'>b$ for all $j<i$ and
\begin{equation}\label{ineq'}
|y_i'-i|<1/10i, \mbox{ where $l_i'Q=(x_i',y_i')$.}
\end{equation}
The set $U'$ exists by Vinogradov's Theorem again. Define the set~$T_1$ to
consist of all prime divisors of
the terms~$B_{l_i'}, i\in\mathbb N$. The set~$T_1$ contains only good
reduction primes, also $T_1\cap S_1=\emptyset$
and $T_1\cap S_2 \neq \emptyset$. Now choose~$T_2$ as follows:
\begin{equation}
T_2=\{p_l':l_i'\neq l \mbox{ prime} \}\cup\{p_{l_i'l_j'}:1\leq j\leq i\}\cup \{p_{ll_i'}:l\in L, i\in \mathbb N \}.
\end{equation}
Then $T_2$ is disjoint from $T_1 \cup S_2$ but it has non-empty
intersection with $S_1$. Now let~$S$ denote any recursive set containing~$S_1\cup T_2$
but disjoint from $S_2\cup T_1$, for example,~$S=S_1\cup T_2$. Then~$S$ will
contain~$S_1$ and be disjoint from~$S_2$. Let~$T$ be the complement of $S$.
The set~$T$ will necessarily contain~$S_2\cup T_1$ and be disjoint from~$S_1\cup T_2$.
It follows that~$T$ will contain~$T_1$ and be disjoint from~$T_2$. A Venn diagram
helps to explain the relationship between these sets.
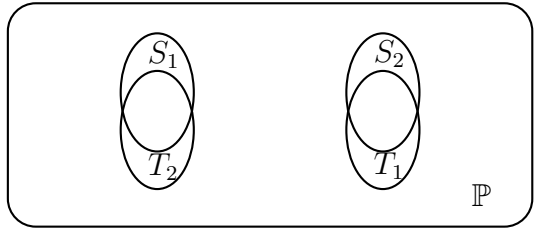
\begin{center}
\begin{pspicture}(-1,0.5)(6,4.5)
\psframe[framearc=0.25](0,1)(7,4)
\psellipse(2,2.3)(.5,0.8)
\uput[1](1.7,3.3){$S_1$}
\psellipse(2,2.8)(.5,0.8)
\uput[2](1.7,1.8){$T_2$}
\psellipse(5,2.3)(.5,0.8)
\uput[3](4.7,3.3){$S_2$}
\psellipse(5,2.8)(.5,0.8)
\uput[4](4.7,1.8){$T_1$}
\uput[5](6,1.4){$\mathbb P$}
\end{pspicture}
\end{center}
The undecidability results follow exactly as before and this completes the proof of Theorem~\ref{h10}.

\subsection{Recursive sets}
The sets
of primes~$S_1,S_2,T_1,T_2$  contain only good reduction primes.
The sets~$U$ and~$U'$ are recursive because the members form a
strictly increasing sequence, the terms of which can be computed in order.
In what follows, let $p>0$ denote a prime of good reduction, and let~$n_p$ denote the
order of~$Q$ mod~$p$. Now~$p|B_{l_i}$ for some~$i$ if and only if $n_p\in U$, which can be checked
because~$U$ is recursive. It follows that~$S_1$ is recursive.
To see if $S_2$ is recursive, first show how to check if $p=p_l$ for some $l \notin U$.
Factorizing $E_p=|E(\mathbb F_p)|$,
one can decide if there is
a prime factor $l|E_p$ such that $l\notin U$ because~$U$ is recursive,
then check if $p=p_l$.
To see if $p=p_{l_il_j}$ for some $1\leq j\leq i$,
factorize~$n_p$ to see if it is the product of two elements~$l_i,l_j \in U$ and~$p=p_{l_il_j}$.
The latter condition can be checked by factorizing earlier terms: in fact, only~$B_{l_i}$
and~$B_{l_j}$ need to be checked. Checking to see if $p=p_{ll_i}$ for some $l\in L$ is similar.
The set~$L$ is recursive because membership can be determined as follows: $l\in L$ if and only if~$B_l=1$.
This completes the proof that~$S_2$ is recursive. The proofs for~$T_1$ and $T_2$ are almost identical, except
that one checks for the second largest prime factor.
\end{proof}

\end{document}